\documentclass[11pt]{article}
\usepackage{latexsym}
\newcounter{example} 
\newtheorem{theorem}{Theorem}  
\newtheorem{lemma}[theorem]{Lemma}

\newtheorem{definition}[theorem]{Definition}

\newcommand{\N}{{\bf N}}

\newcommand{\D}{\mathcal{D}}
\newcommand{\F}{\mathcal{F}}
\renewcommand{\H}{\mathcal{H}}
\newcommand{\g}{\gamma}
\newcommand{\proof}{\noindent \textbf{Proof} \,}
\newcommand{\qed}{\quad \hspace*{\fill} $\Box$ \par} 

\newcommand{\fq}{\mathbf{F}_{q}}

\begin{document}

\title{
   \textbf{\textit{q}-Exponential Families}
}
\author{
    Kent E. Morrison \\
    Department of Mathematics \\
    California Polytechnic State University \\
    San Luis Obispo, CA 93407 \\
    \\
     \texttt{kmorriso@calpoly.edu 
}
   }
    
\date{12 February 2004}   

\maketitle 

\begin{abstract} We develop an analog of the exponential
families of Wilf in which the label sets are finite dimensional vector
spaces over a finite field rather than finite sets of positive
integers.  The essential features of exponential families are
preserved, including the exponential formula relating the deck
enumerator and the hand enumerator.    
\end{abstract}

\section{Introduction and Definitions}
\noindent In this paper we analogize Wilf's approach to labelled
counting in \cite{Wilf95} based on exponential families.  Notation and
definitions follow his as closely as possible.  This work is further
elaboration of the subset-subspace analogy that has long been a rich
source in enumerative combinatorics.  See Kung \cite{Kung95} for an
historical survey.

Let $\fq$ the finite field of order $q$ and ${\fq}^{(\N)}$ the vector
space of countable dimension over $\fq$ whose elements are infinite
sequences $(a_{1},a_{2},\ldots)$ with a finite number of non-zero
entries.  Let $e_{1},e_{2},\ldots$ be the standard basis and define
$E_{n}$ to be the span of $e_{1},\ldots,e_{n}$.  Let $P$ be an
abstract set of `pictures.'

\begin{definition} 
   A \textbf{q-card} $Q(V,p)$ is a pair consisting of a subspace 
   $V \subset {\fq}^{(\N)}$ and a picture $p \in P$. $V$ is the
   \textbf{label space}
   space and the \textbf{dimension} of $Q$ is $\dim V$. A card is 
   \textbf{standard} if its label space is $E_{n}$.
\end{definition}

\begin{definition}  
   A \textbf{q-hand} is a finite set of $q$-cards whose label spaces 
   form a direct sum decomposition of $E_{n}$ for some $n$. The 
   \textbf{dimension} of the hand is $n$.
\end{definition}

\begin{definition} 
   A $q$-card $Q(V',p)$ is a \textbf {relabeling} of $Q(V,p)$ if $\dim 
   V'=\dim V$. (Thus, the pictures must be the same.) The 
   \textbf{standard relabeling} of $Q(V,p)$ is $Q(E_{\dim V},p)$.
\end{definition}

\begin{definition} 
   A \textbf{q-deck} $\mathcal{D}$ is a finite set of standard $q$-cards whose 
   dimensions are the same and whose pictures are different. The 
   dimension of the deck $\dim \mathcal{D}$ is the common dimension 
   of the cards.
\end{definition}

Generally, we will omit the \emph{q}-prefix and refer to ``cards, 
hands, and decks'' when the context makes it clear.
\begin{definition} 
   A \textbf{q-exponential family} $\F$ is a collection of $q$-decks 
   $\D_{1},\D_{2},\ldots$ where $\D_{n}$ is a $q$-deck of dimension 
   $n$, possibly empty.
\end{definition}

\begin{definition}
   Given two $q$-exponential families $\F'$ and $\F''$ the 
   \textbf{merger} is the $q$-exponential family $\F' \oplus \F''$ 
   whose $q$-deck of dimension $n$ is the disjoint union of the 
   $q$-decks $\D'_{n}$ and $\D''_{n}$.
\end{definition}

For notational simplicity we define $\gamma_{n}=|\mathrm{GL}_{n}(\fq)|$ to be the order of 
the general linear group, which is given by
\[ \gamma_{n} =\prod_{0\leq i \leq n-1} (q^{n}-q^{i}) .\]
Define the \emph{hand enumerator} 
\[ \H(x,y) = \sum_{n,k}h(n,k)\frac{x^{n}}{\gamma_{n}}y^{k}  \]
and the \emph{deck enumerator}
\[  \D(x) = \sum_{n}\frac{d_{n}}{\gamma_{n}}x^{n} , \]
where $h(n,k)$ is the number of hands of dimension $n$ with $k$ cards 
and $d_{n}$ is the number of cards in $\D_{n}$. We define $h(0,0)$ to 
be 1 and will see later why that is necessary. However, 
$h(n,0)=h(0,k)=0$ for $n > 0$ and $k > 0$.

\section{Counting}
The one variable hand enumerator 
\[ \H(x) = \H(x,1)= \sum_{n}h(n) \frac{x^{n}}{\g_{n}} \]
counts the hands without regard to the number of cards. Here, 
$h(n)=\sum_{k}h(n,k)$ is the number of hands of dimension $n$.
The fundamental lemma of $q$-labeled counting is key in deriving the 
analog of the exponential formula.
\begin{lemma}
   Let $ \F'$ and $\F''$ be $q$-exponential families and $\F = \F' 
   \oplus \F''$ the merger of the two. Then the two-variable hand 
   enumerators are related by 
   \[ \H(x,y) = \H'(x,y)\H''(x,y) .\]
\end{lemma}
\proof  A hand of dimension $n$ with $k$ cards has label spaces
that give a splitting of $E_{n}$.  Those cards that come from decks in
$\F'$, say $k'$ of them have label spaces that form a subspace $V'$ of
dimension $n'$, and those that come from decks of $\F''$, the
remaining $k-k''$ cards, have label spaces that form a complementary
subspace $V''$ of dimension $n-n'$. Thus, we select a hand from $\F'$ 
and a hand from $\F''$ and relabel the cards by selecting a splitting 
$E_{n} = V' \oplus V''$. The number of ways to select $V'$ and $V''$ 
of dimensions $n'$ and $n-n'$ is $\frac{\g_{n}}{\g_{n'}\g_{n-n'}}$, 
because $\mathrm{GL}_{n}(\fq)$ acts transitively on the set of pairs 
$(V',V'')$ such that $E_n = V' \oplus V''$ and $\dim V' = n'$, $\dim 
V''=n-n''$. The stabilizer subgroup is isomorphic to 
$\mathrm{GL}_{n'}(\fq) \times \mathrm{GL}_{n-n'}(\fq)$. Therefore 
\[ 
   h(n,k) = \sum_{n',k'}
   \frac{\g_{n}}{\g_{n'}\g_{n-n'}}h'(n',k')h''(n-n',k-k') .
\]
Here we  need $h(0,0)=h'(0,0)=h''(0,0)=1$ to allow for the 
possibility of making up a hand in the merger by choosing no cards 
from one of the families. We might also point out that $\g_{0}=1$ 
because the zero map on the vector space of dimension zero is 
invertible.

In the product 
\[ \H'(x,y)\H''(x,y)=\sum_{n',k'}h'(n',k')\frac{x^{n'}}{\g_{n'}}y^{k'}
          \sum_{n'',k''}h''(n'',k'')\frac{x^{n''}}{\g_{n''}}y^{k''}  \]
 the coefficient of $x^{n}y^{k}$ is
\[ \sum_{n'+n''=n,k'+k''=k} 
h'(n',k')h''(n'',k'')\frac{1}{\g_{n'}\g_{n''}} . \]
Substituting $n''=n-n'$ and $k''=k-k'$ shows that this is exactly 
$h(n,k)/\g_{n}$.  \qed

We continue on to prove the exponential formula using the three steps 
that Wilf calls \emph{the trickle, the flow, and the flood}.

For the first step, we consider a family consisting of one deck 
$\D_{r}$ with one card. Thus, $d_{r}=1$ and all other $d_{i}=0$. A 
hand with $k$ cards has dimension $n=kr$. The label spaces for a hand 
form a splitting of $E_{n}$ into $k$ subspaces each of dimension $r$. 
Thus the number of hands of dimension $n$ is the number of such splittings
$E_{n}=V_{1} \oplus \cdots \oplus V_{k}$. Changing the order of the 
factors does not change the splitting, so we can count the number of 
ordered splittings and divide by $k!$. The group 
$\mathrm{GL}_{n}(\fq)$ acts transitively on the set of ordered 
splittings with stabilizer subgroup isomorphic to 
$\mathrm{GL}_{r}(\fq)^{k}$, and so the number of $n$-dimensional 
hands (where $n=kr$) is 
\[ h(kr,k)=\frac{1}{k!}\frac{\g_{kr}}{\g_{r}^{k}} . \]
The hand enumerator is 
\begin{eqnarray*}
    \H(x,y) & = & \sum_{k} \frac{1}{k!}\frac{\g_{kr}}{\g_{r}^{k}} 
    \frac{x^{kr}}{\g_{kr}}y^{k}  \\
        & =&  \sum_{k}\frac{1}{k!}\left(\frac{yx^{r}}{\g_{r}}\right)^{k}\\
	   &=& \exp\left(\frac{yx^{r}}{\g_{r}}\right) \\
	   &=& \exp(y \D(x) ) .
\end{eqnarray*}

The second step handles the case of a single non-empty deck $\D_{r}$ 
with $d_{r}$ cards. If we denote the family in step one by $\F_{1}$, 
then we are now considering the repeated merger
$\F_{1} \oplus\cdots \oplus \F_{1}$, 
in which there are $d_{r}$ factors. Therefore, the hand 
enumerator is the $d_{r}$-th power of the hand enumerator in step one, 
namely 
\[ \H(x,y)= \left(\exp\left(\frac{yx^{r}}{\g_{r}}\right)\right)^{d_{r}} = 
\exp\left(\frac{yd_{r}x^{r}}{\g_{r}}\right) . \]
The deck enumerator is $\D(x) = \frac{d_{r}x^{r}}{\g_{r}}$, and so we 
have $\H(x,y)=\exp(y\D(x))$.

In the third step we look at a full $q$-exponential family as a merger 
of $\F_{r}$, $r=1,2,\ldots$, where $\F_{r}$ is the family consisting 
of just one non-empty deck $\D_{r}$. Then the hand enumerator of the 
infinite merger is the infinite product of the hand enumerators
\begin{eqnarray*}
   \H(x,y) &=& \prod_{r} \exp\left(\frac{yd_{r}x^{r}}{\g_{r}}\right)\\
     &=& \exp\left(\sum_{r}\frac{yd_{r}x^{r}}{\g_{r}} \right) .
\end{eqnarray*}
The sum in the last line is the deck enumerator of the family, and so 
we have proved
\begin{theorem}[Analog of the exponential formula]
   Let $\F$ be a $q$-exponential family with hand enumerator 
   $\H(x,y)$ and deck enumerator $\D(x)$. Then
   \[ \H(x,y)= \exp(y\D(x)). \] 
\end{theorem}

In \cite[\S 3.13]{Wilf95} Wilf defines the polynomials $\phi_{n}(y) = 
\sum_{k}h(n,k)y^{k}$. The exponential formula can be written
\[ e^{y\D (x)} = \sum_{n \geq 0}\frac{\phi_{n}(y)}{n!}x^{n} , \]
from which it follows that 
it follows that these polynomials satisfy 
\[ \phi_{n}(u+v)=\sum_{m}{n \choose m}\phi_{m}(u)\phi_{n-m}(v)  \]
so that they are \emph{polynomials of binomial 
type} as defined by Mullin and Rota \cite{MullinRota70}.

For $q$-exponential families there is an analog: for the polynomials 
$\phi_{n}(y)$ defined as above, the $q$-exponential formula just 
proved can be written
\[ e^{y\D(x)} = \sum_{n /geq 0}\frac{\phi_{n}(y)}{\g_{n}}x^{n} . \]
It follows that 
\[ \phi_{n}(u+v) = 
   \sum_{m}\frac{\g_{n}}{\g_{m}\g_{n-m}}\phi_{m}(u)\phi_{n-m}(v) . \]
   
\section{\textit{q}-Analog of the Stirling subset numbers}
The Stirling subset number ${n \brace k}$, or Stirling number of 
the second kind, counts the number of partitions of a set of size $n$ 
into $k$ non-empty subsets. The direct sum of subspaces is analogous 
to the disjoint union of subsets, and so we define ${n \brace k}_{q}$ 
to be the number of splittings of $E_{n}$ into $k$ non-zero 
subspaces. We regard a splitting as a hand of dimension $n$ with $k$ 
cards in the $q$-exponential family $\F$ with decks $\D_{r}$ 
consisting of just one card each. (The picture set is irrelevant. We 
can assume there is just one picture and that it is the same on each 
card.) Then $h(n,k)= {n \brace k}_{q}$ and 
\begin{eqnarray*}
\H(x,y) &=& \sum_{n,k}{n \brace k}_{q} \frac{x^{n}}{\g_{n}}y^{k} \\
\D(x) &= & \sum_{r}\frac{x^{r}}{\g_{r}} .
\end{eqnarray*}
From the exponential formula 
\begin{eqnarray*}
   {n \brace k}_{q} & = & \g_{n} [x^{n}y^{k}]\H(x,y) \\
   &=& \g_{n} [x^{n}y^{k}]\sum_{j}\frac{y^{j}\D(x)^{j}}{j!}\\
   &=& \g_{n} [x^{n}]\frac{\D(x)^{k}}{k!} \\
   &=& \g_{n}[x^{n}] 
   \frac{1}{k!}\left(\sum_{r}\frac{x^{r}}{\g_{r}}\right)^{k} \\
   &=& \g_{n}\frac{1}{k!} \sum_{\stackrel{n_{1}+\cdots+n_{k}=n}{n_{i}\geq 1}}
        \frac{1}{\g_{n_{1}}\cdots\g_{n_{k}}} .
\end{eqnarray*} 
Thus we see that 
\[ 
   {n \brace k}_{q} = \frac{1}{k!} \sum_{\stackrel{n_{1}+\cdots+n_{k}=n}{n_{i}\geq 1}}
        \frac{\g_{n}}{\g_{n_{1}}\cdots\g_{n_{k}}} .
\]
Compare this with the exponential formula for the Stirling subset numbers
\[  \sum_{n,k}{n \brace k}\frac{x^{n}}{n!}y^{k}=\exp(y(e^{x}-1))  \] 
to see that we have a perfect analog of the formula
\[  {n \brace k} = \frac{1}{k!} \sum_{\stackrel{n_{1}+\cdots+n_{k}=n}{n_{i}\geq 1}}
        \frac{n!}{{n_{1}!}\cdots {n_{k}!}}  . \]

The Bell number $b(n)=\sum_{k}{n \brace k}$ counts the number of set
partitions of a set of size $n$ and therefore has a $q$-analog
$b_{q}(n)=\sum_{k}{n \brace k}_{q}$, which counts the number of direct
sum decompositions of an $n$-dimensional vector space over $\fq$. 
These are the coefficients in the one-variable hand enumerator
$\H(x)=\H(x,1)=\sum_{n}b_{q}(n)\frac{x^{n}}{\g_{n}}$.  By setting
$y=1$ in the $q$-exponential formula we get 
\[
\sum_{n}b_{q}(n)\frac{x^{n}}{\g_{n}}=\exp\left(\sum_{r \geq
1}\frac{x^{r}}{\g_{r}} \right) . 
\] 
This formula for the exponential
generating function of the $q$-Bell numbers first appears as Example
11 in \cite{BenderGoldman71}.  It can also be derived from the
exponential formula in \cite{Stanley78}.
(Look at Example 2.2 and then use Corollary 3.3 with $f(n)=1$ for $n > 0
$. Notice that $n!M(n)=\g_{n}/(q-1)^{n}$ and substitute $x/(q-1)$ for
$x$.)

\section{Diagonalizations}

Consider the $q$-exponential family with just one non-empty deck
$\D_{1}$ with $q$ cards.  A card in $\D_{1}$ has the label $E_{1}$ and
a picture $\alpha$, which is an element of $\fq$.  Then a hand of
dimension $n$ must consist of $n$ cards whose label spaces form a
splitting of $E_{n}=V_{1} \oplus \cdots \oplus V_{n}$ into
one-dimensional subspaces and whose pictures
$\alpha_{1},\ldots,\alpha_{n}$ represent a diagonal matrix with
respect to that splitting.  Thus, each hand of dimension $n$
represents a diagonalization of a diagonalizable $n \times n$ matrix.

The deck enumerator is \[ \D(x) = \frac{q}{\gamma_{1}}x =
      \frac{q}{q-1}x ,\]
and the one variable hand enumerator is 
\[ \H(x)=\exp\left( \frac{q}{q-1}x\right) .\]
Since every hand of dimension $n$ has $n$ cards, there is no more 
information in the two variable hand enumerator. Then, $h(n)$, the 
number of diagonalizations of $n \times n$ matrices is
\[ h(n)= \frac{\gamma_{n}}{n!}\left(\frac{q}{q-1} \right)^{n} .\]

A variation of this family changes the picture set slightly so that 
only non-zero $\alpha$ in $\fq$ are used. In this case a hand 
corresponds to a diagonalization of an invertible matrix. Since 
$d_{1}=q-1$, the deck and hand enumerators become $\D(x)=x$ and 
$\H(x)=e^{x}$. Thus, the number of diagonalizations of invertible $n 
\times n$ matrices is \[ h(n)= \frac{\gamma_{n}}{n!} . \]
A direct proof without generating functions can be given for this 
result, but we leave that as an exercise for the reader.

In order to count diagonalizable matrices rather than 
diagonalizations we have to go beyond exponential families. The 
problem can be seen in the simplest of cases. The $2 \times 2$ 
identity matrix has as many diagonalizations as there are splittings 
$E_{2}=V_{1}\oplus V_{2}$ into two one-dimensional subspaces. 
However, a $2 \times 2$ matrix with distinct eigenvalues has only one 
diagonalization. We take this up in the next section.

\section{Beyond exponential families}

\subsection{Diagonalizable matrices}
An $n \times n$ diagonalizable matrix gives a unique decomposition of
$E_{n}$ into a direct sum of eigenspaces, one for each eigenvalue. 
However, in general it does not give a unique decomposition of $E_{n}$
into a sum of one-dimensional spaces.  This fact forces us to go
beyond $q$-exponential families.  We would like a hand of dimension
$n$ to correspond to the eigenspace decomposition of a diagonalizable
$n \times n$ matrix.  Thus, a card needs a label space $V \subset
E_{n}$ and an eigenvalue `picture' $\alpha \in \fq$.  However, a
picture cannot be used more than once and a hand cannot contain more
than $q$ cards.  In spite of this difficulty, we can see the family of
diagonalizable matrices as the merger of $q$ families, one for each
possible eigenvalue.  The fundamental counting lemma still holds for
the merger of families and this will allow us to count.

\newcommand{\diag}{\mathrm{diag}}
For each $\alpha$ in $\fq$ we define the family whose hands consist of 
single cards with label space $E_{n}$ and picture $\alpha$, $n \geq 
0$. We define the hand enumerator $\H_{\alpha}(x,y)$ to be the 
generating function
\[ \H_{\alpha}(x,y)= 1 + y\sum_{n \geq 1} \frac{x^{n}}{\gamma_{n}} . \]
Notice that $y$ appears only with exponent 1 because every hand has 
only one card. Now we can form the product of these hand enumerators 
as $\alpha$ ranges over $\fq$, and by the fundamental of $q$-labeled 
counting, that product will be the hand enumerator for the family of 
diagonalizable matrices:
\[ \H_{\mathrm{diag}}(x,y)= \left(1 + y\sum_{n \geq 1} 
\frac{x^{n}}{\gamma_{n}} \right)^{q}  . \]
Here the coefficient of $x^{n}y^{k}$, when multiplied by $\gamma_{n}$ 
is the number of $n \times n$ diagonalizable matrices with $k$ 
distinct eigenvalues. Setting $y=1$ gives the one variable enumerator 
$\H_{\diag}(x)= \sum_{n \geq 0}\frac{h_{\diag}(n)}{\gamma_{n}}x^{n}$. 
Hence, 
\begin{eqnarray*}
 h_{\diag}(n) &=& \gamma_{n} [x^{n}] \left( \sum_{m \geq 
0}\frac{x^{m}}{\gamma_{m}} \right)^{q}  \\
&=& \sum_{n_{1}+\cdots + n_{q}=n}\frac{\gamma_{n}}{\gamma_{n_{1}}\cdots 
\gamma_{n_{q}}} . 
\end{eqnarray*}

\subsection{Projections}
A projection is a matrix $P$ such that $P^{2}=P$.  A projection is
completely determined by the pair of complementary subspaces, the
kernel of $P$ and the image of $P$.  Thus, the number of $n \times n$
projections can be obtained from the $q$-Stirling number 
${n \brace 2}_{q}$, which counts the number of decompositions into two
non-trivial subspaces.  First we have to multiply that by two because
the order of the two subspaces now matters and then we have to add two
to count the trivial decompositions corresponding to the identity
matrix and to the zero matrix.  Thus, the number of $n \times n$
projections is $ 2 {n \brace 2}_{q} +2$.  Recalling the formula from
section 3 for the $q$-Stirling numbers we see that the number of
projections is
\begin{eqnarray*}
2{n \brace 2}_{q}+2 &=& \sum_{\stackrel{n_{1}+n_{2}=n}{n_{i}\geq 1}}
        \frac{\g_{n}}{\g_{n_{1}} \g_{n_{2}}}  + 2  \\
	&=&  \sum_{\stackrel{n_{1}+n_{2}=n}{n_{i}\geq 0}}
        \frac{\g_{n}}{\g_{n_{1}} \g_{n_{2}}} .
\end{eqnarray*}
The 2 on the right is incorporated into the sum by allowing $n_{1}$ 
and $n_{2}$ to be 0.

An alternative approach is to use a variation of the family of
diagonalizable matrices.  A projection is a diagonalizable matrix
whose only eigenvalues are 0 or 1.  The family of projections is the
merger of the two families with $\alpha=0$ and $\alpha = 1$.  Then the
one variable enumerator for the family of projections is
\[  \H_{\mathrm{pr}}(x)= 
    \left( \sum_{m \geq 0}\frac{x^{m}}{\g_{m}} \right)^{2} ,
\]
and the number of projections within the $n \times n$ matrices is
\[  h_{\mathrm{pr}}(n)= \sum_{j=0}^{n} \frac{\g_{n}}{\g_{j}\g_{n-j}} .
\]

\subsection{A \textit{q}-analog of the Stirling cycle numbers}
The Stirling cycle number ${n \brack k}$, or Stirling
number of the first kind, is the number of permutations on $n$ letters
having $k$ cycles.  Replacing a set of size $n$ by a vector space of
dimension $n$ over $\fq$ and the permutation group $S_{n}$ by the
general linear group $\mathrm{GL}_{n}(\fq)$, a reasonable $q$-analog of
the Stirling cycle numbers defines ${n \brack k}_{q}$ to be the number
of invertible $n$ by $n$ matrices that decompose into the direct sum
of $k$ cyclic summands.  In order to make this definition precise we
specify that we are counting the number of summands in the primary
rational canonical form.  Recall, that each endomorphism of a finite
dimensional vector space decomposes uniquely into cyclic summands
whose matrix representations are those of companion matrices of powers
of irreducible polynomials.  (Without that specification we have the
ambiguity, for example, of a $2 \times 2$ matrix with distinct
eigenvalues being itself cyclic and also decomposing as the direct sum of
two $1 \times 1$ matrices.  In this case the number of primary cyclic
summands is two.)

We construct the appropriate family as the infinite merger of families
$\F_{\phi}$, where $\phi$ is an irreducible monic polynomial over
$\fq$.  A hand in $\F_{\phi}$ is a matrix whose characteristic
polynomial is a power of $\phi$.  The conjugacy class of such a matrix
is determined by the multiplicity data $(b_{i})$, where $b_{i}$ is the
number of copies of the companion matrix $C(\phi^{i})$ occurring.  The
dimension of the hand is $\sum_{i}ib_{i}(\deg \phi)$, because
$C(\phi^{i})$ is a matrix of dimension $i(\deg \phi)$.  Note that the
multiplicity data is a partition of the integer $\sum_{i}ib_{i}$.  The
number of parts in the partition $b$ is the sum $k=\sum_{i}b_{i}$, and
this is the number of cyclic summands in the rational primary
canonical form.

To describe the hand enumerator for $\F_{\phi}$ we need the number of
matrices whose conjugacy class is determined by the multiplicity data
$b$.  The cardinality of this conjugacy class is
$\g_{n}/c_\phi(b)$, where $n=\sum_{i}ib_{i}(\deg \phi)$ and
$c_\phi(b)$ is the order of the stabilizer subgroup under the
conjugation action of the rational canonical form associated to 
$\phi$ and $b$.   Summing over all
partitions $b$ we get the hand enumerator
\[ 
  \H_{\phi}(x,y) = 
  \sum_{b}\frac{x^{\sum ib_{i}\deg \phi}y^{\sum b_{i}}}{c_\phi(b)} . 
\]
Note that the exponent in the $x$ variable is the dimension and the 
exponent in the $y$ variable is the number of cyclic summands. The 
proof of the
formula for $c_{\phi}(b)$ can be found in \cite{Kung81}. We give the 
formula here for completeness, although we make no further use of it. 
For the partition $b=(b_{1},b_{2},\ldots)$ define 
\[ d_{i}= b_{1}+2b_{2}+\cdots + (i-1)b_{i-1} + i(b_{i}+b_{i+1} + \cdots ) . \]
In the Ferrers diagram of $b$ consisting or $b_{i}$ rows with $i$ 
boxes, $d_{i}$ is the number of boxes in the first $i$ columns. Then
\[ c_{\phi}(b) = \prod_{i}\prod_{j=1}^{b_{i}}
                 (q^{d_{i}\deg \phi}-q^{(d_{i}-j)\deg \phi}) .
\]

Forming the merger of the families $\F_{\phi}$ as $\phi$ ranges over
the irreducible monic polynomials in $\fq[t]$ results in the family of
all square matrices over $\fq$ decomposed into their rational primary
components.  In order to have only the invertible matrices, we simply
omit the polynomial $\phi(t)=t$ in forming the merger, because a
matrix is invertible if and only if its characteristic polynomial is
not divisible by $t$.  With this family we have a factorization of the
generating function for the $q$-analog of the Stirling cycle numbers
\[
   \sum {n \brack k}_{q}\frac{x^{n}}{\g_{n}}y^{k}
   =\prod_{\phi \neq t}\left(1+
     \sum_{b}\frac{x^{\sum ib_{i}\deg \phi}y^{\sum b_{i}}}{c_\phi(b)} 
     \right)
\]
Note the constant 1 in each factor to allow for $\phi$ not occurring 
in the primary decomposition of a matrix. This dictates that ${0 
\brack 0}_{q}$ should be 1 just as it is for the ordinary Stirling 
numbers.

This factorization of the generating function for these $q$-Stirling 
numbers was first obtained by Kung \cite[p. 148]{Kung81}. (To compare 
the formulas replace the $x$ and $y$ in this paper with $u$ and $x$ in 
Kung's paper, break the sum over all partitions into a sum over $j 
\geq 1$ and an inner sum over partitions of $j$, and group the factors 
for all irreducible polynomials of the same degree.) In that 
important paper he defined the cycle index for groups of 
automorphisms of finite dimensional vector spaces over finite fields 
and found its basic properties. The ordinary generating function for 
the cycle index of the full general linear group has a factorization 
\cite[Lemma 1]{Kung81} that specializes to the factorization above 
for the $q$-Stirling cycle numbers. Extensive use of the vector 
space cycle index is made in later papers of Stong \cite{Stong88} and 
Fulman \cite{Fulman02} to study the asymptotic combinatorics of the 
canonical form as the matrix size goes to infinity.
Some of this work is also summarized in \cite{Morrison99}.

\end{document}